%
\nopagenumbers
\input amstex
\documentstyle {amsppt}
\magnification=\magstep1
\parindent=20pt
\hfuzz=10pt

\centerline {\bf A division algorithm for the free left distributive
algebra}
\vskip 1em
\centerline {Richard Laver\footnote{Supported by NSF Grant DMS
9102703.}}
\centerline {University of Colorado, Boulder}
\vskip 2em
\baselineskip=20pt

In this paper we extend the normal form theorem, for the free algebra
${\Cal A}$ on one generator $x$ satisfying the left distributive law
$a(bc) = (ab)(ac)$, which was proved in [5].  As part of the proof that
an algebra of elementary embeddings from set theory is isomorphic to
${\Cal A}$, facts about ${\Cal A}$ itself were established.  Theorem 1
summarizes some known facts about ${\Cal A}$, including P. Dehornoy's
independent work on the subject.  After that the main
theorem, about putting members of ${\Cal A}$ into ``division form,''
will be proved with the help of versions of lemmas of [5] and one of
the normal forms of [5]. 
  
Let $\cdot$ denote the operation of ${\Cal A}$.  These forms take place
not in ${\Cal A}$ but in a larger algebra ${\Cal P}$ which involves
additionally a composition operation $\circ$.  Let $\sum$ be the set of
laws $\{ a\circ (b\circ c) = (a\circ b)\circ c$, $(a\circ b)c = a(bc)$,
$a(b\circ c) = ab \circ ac$, $a\circ b = ab \circ a\}$.  ${\Cal P}$ is
the free algebra on the generator $x$ satisfying $\sum$.  $\sum$ implies
the left distributive law, and $\sum$ is a conservative extension of the
left-distributive law (if two terms in the language of ${\Cal A}$ can be
proved equal using  $\sum$, then they can be proved equal using just
the left-distributive law).  So we may identify ${\Cal A}$ as a subalgebra
of ${\Cal P}$ restricted to $\cdot$.   If $p_0 ,p_1 ,\ldots ,p_n \in 
{\Cal P}$, write
$p_0 p_1 \cdots p_n$ (respectively, $p_0 p_1 \cdots p_{n-1} \circ p_n$)
for $(((p_0 p_1 )p_2 )\cdots p_{n-1} )p_n$ (respectively, $(((p_0 p_1
)p_2 )\cdots p_{n-1} )\circ p_n$).  Write $w = p_0 p_1 \cdots
p_{n-1} \ast p_n$ to mean that either $w = p_0 p_1 \cdots p_n$ or $w =
p_0 p_1 \cdots p_{n-1}  \circ p_n$.  Make these conventions also for
other algebras on operations $\cdot$ and $\circ$.

For $p \in {\Cal P}$ let $p^1 =p$, $p^{n+1} =p \circ p^n$; let $p^{(0)}
= p$, $p^{(n+1)} = pp^{(n)}$.  Then $p^{(n+1)} = p^{(i)} p^{(n)}$ for
all $i \leq n$, by induction using the left-distributive law.

For $p, q \in {\Cal P}$ let $p < q$ if $q$ can be written as a term of
length greater than one in the operations $\cdot$ and $\circ$,involving
members of ${\Cal P}$ at least one of which is $p$.  
Write $p <_L q$ if
$p$ occurs on the left of such a product:  ~$q = pa_0 a_1 \cdots
a_{n-1}\ast a_n$ for some $n \geq 0$.  Then $<_L$ and $<$ are
transitive.  If $a, b \in {\Cal A}$ and $a <_L b$ in the sense of ${\Cal
P}$, then $a <_L b$ in the sense of ${\Cal A}$; and similarly for $<$.

In [5] it was shown,via the existence of normal forms for the members of        ${\Cal P}$,that $<_L$ linearly orders ${\Cal P}$ and ${\Cal A}$.
The proof of
part of that theorem, that $<_L$ is irreflexive, used a large cardinal
axiom (the  existence, for each $n$, of an $n$-huge cardinal).  Dehornoy
([1], [2]) by a different method independently proved in ZFC that for all $a, b \in {\Cal A}$
at least one of $a <_L b$, $a = b$, $b <_L a$ holds. 
Recently ([3]) he has found a proof of the irreflexivity of $<-L$ in ZFC.
Dehornoy's theorem has the consequence that facts about ${\Cal P}$ 
(Theorem 1 below(parts (v)--(viii)),and the normal and division forms
in [5] and this paper)
which have previously been 
known from a large cardinal assumption(that is,from irreflexivity),are
provable in ZFC.

For $u, v$ terms in the language of $\cdot$ in the variable $x$, let $u
\rightarrow v$ ([1]) mean that $u$ can be transformed into $v$ by a finite
number of substitutions, each consisting of replacing a term of the form
$a(bc)$ by $(ab)(ac)$.

For $\lambda$ a limit ordinal, let ${\Cal E}_\lambda$ be the set of
elementary embeddings $j$; $(V_\lambda ,\epsilon ) \rightarrow (V_\lambda
,\epsilon )$, $j$ not the identity.  For $j, k \in {\Cal E}_\lambda$,
let $jk = \underset \alpha < \lambda \to\bigcup j(k \cap V_\alpha )$
and let $j \circ k$ be the composition of $j$ and $k$.  Then the
existence of a $\lambda$ such that ${\Cal E}_\lambda \not= \phi$ is a
large cardinal axiom, if $j, k \in {\Cal E}_\lambda$, then $jk$, $j\circ k
\in {\Cal E}_\lambda$, and $({\Cal E}_\lambda ,\cdot ,\circ )$ satisfies
$\sum$.  For $j \in {\Cal E}_\lambda$ let ${\Cal A}_j$ be the closure of
$\{ j\}$ under $\cdot$ and let ${\Cal P}_j$ be the closure of $\{ j\}$
under $\cdot$ and $\circ$.

Some facts relating ${\Cal P}$ to ${\Cal A}$, such as the conservativeness
of $\sum$ over the left-distributive law, may be found in [5].

\proclaim{Theorem 1}  (i) If $r <_L s$, then $pr <_L p\circ r <_L ps$.
\item {(ii)} $x \leq_L p$ for all $p \in {\Cal P}$, $<_L$ is not well
founded.
\item {(iii)} For all $p, q \in {\Cal P}$ there is an $n$ with $p^{(n)}
> q$.
\item {(iv)} The rewriting rules for ${\Cal A}$ are confluent, i.e., if
$u, v$ are terms in the language of $\cdot$ in the variable $x$, and $u
\equiv v$ via the left distributive law, then for some $w$, $u \rightarrow
w$ and $v \rightarrow w$.
\item {(v)} $<_L$ is a linear ordering of ${\Cal A} , {\Cal P}$.
\item {(vi)} For $p, q, r \in {\Cal P}$, $pq = pr \Leftrightarrow q = r$,
$pq <_L pr \Leftrightarrow q <_L r$.
\item {(vii)} The word problems for ${\Cal A}$ and ${\Cal P}$ are
decidable.
\item {(viii)} $<_L = <$ on ${\Cal A} , {\Cal P}$.
\item {(ix)} For no $k_0 ,k_1 ,\ldots ,k_n \in {\Cal E}_\lambda$ ($n >
0$) is $k_0 = k_0 k_1 \cdots k_{n-1} \ast k_n$.
\item {(x)} For all $j \in {\Cal E}_\lambda$, ${\Cal A}_j \cong {\Cal
A}$, ${\Cal P}_j \cong {\Cal P}$.
\endproclaim

{\bf Remarks.}  (i)--(iii) are quickly proved; for (iii), it may be seen
that 
$p^{(n)} \geq x^{(n)}$ and for sufficiently large $n$, $x^{(n)} \geq q$.
(iv) is Dehornoy's theorem in [2].  
The linear orderings of ${\Cal P}$
and ${\Cal A}$ both have order type $\omega \cdot (1+\eta )$.  
(v) immediately implies (vi) and (vii).  In [5],
(viii) is derived from the normal form theorem; McKenzie derived (viii)
from (v).  (ix) and (x) are proved in [5], (ix)
plus (v) yields (x).

Results connected with critical points of members of ${\Cal A}_j$
appear in [4],[6] and [7].

For $a, b \in {\Cal P}$, let the iterates  $I_n (a,b)$ of $\langle a,b
\rangle$ ($n \geq 1$) be defined by $I_1 (a,b) = a$, $I_2 (a,b) = ab$,
$I_{n+2} (a,b) = I_{n+1} (a,b) I_n (a,b)$.

Call a term $b_0 b_1 \cdots b_{n-1} \ast b_n$, with each $b_i \in {\Cal
P}$, prenormal (with respect to a given ordering $\prec$) if $b_2
\preceq b_0$, $b_3 \preceq b_0 b_1$, $b_4 \preceq b_0 b_1 b_2 ,\ldots 
,b_n \preceq b_0 b_1 \cdots b_{n-2}$, and in the case $\ast = \circ$ and 
$n \geq 2$, $b_n \prec b_0 b_1 \cdots b_{n-1}$.

The main theorem is that for each $p, q \in {\Cal P}$, $q$ can be
expressed in ``$p$-division form,''the natural fact suggested by the
normal forms of [5]. For $p \in {\Cal P}$ the set of
$p$-division form representations of members of ${\Cal P}$, and its
lexicographic linear ordering, are defined as follows.

\proclaim{Lemma 2}  For each $p \in {\Cal P}$ there is a unique set
$p-DF$ of terms in the language of $\cdot$ and $\circ$, in the alphabet
$\{ q \in {\Cal P} : q \leq p \}$, and a linear ordering
$<_{\text{Lex}}$ of $p-DF$, such that
\item {(i)} For each $q \leq_L p$, $q$ (as a term of length one) is in
$p-DF$, and for $q, r \leq_L p$, $q <_{\text{Lex}} r$  if and only if
$q <_L r$. 
\item {(ii)} $w \in p-DF$ iff either $w \leq_L p$, or $w = pa_1 a_1
\cdots a_{n-1} \ast a_n$, where each $a_i \in p-DF$, is prenormal with
respect to $<_{\text{Lex}}$.
\item {(iii)} For $w \in p-DF$ define the associated sequence of $w$
to be $\langle w\rangle$ if $w \leq_L p$, to be $\langle p, a_0 ,a_1
,\ldots ,a_n \rangle$ if $w = pa_0 a_1 \cdots a_n$, and, if $w = pa_0
a_1 \cdots a_{n-1} \circ a_n$, to be (letting $u =  pa_0 a_1 \cdots
a_{n-1}$) $\langle p, a_0 ,a_1 ,\ldots ,a_{n-1} ,a_n ,u, ua_n ,ua_n u,
ua_n u(ua_n ),\ldots \rangle$, that is, the sequence beyond $a_n$ is 
$\langle I_m (u, a_n ) : m \geq 1\rangle$.  Then
for $w, v \in DF$ with associated sequences $\langle w_i : i < \alpha
\rangle$, $\langle v_i : i < \beta \rangle$ ($\alpha ,\beta \leq \omega
\rangle$, $w <_{\text{Lex}} v$ iff either $\langle w_i : i < \alpha
\rangle$ is a proper initial segment of $\langle v_i : i < \beta
\rangle$ or there is a least $i$ with $w_i \not= v_i$, and $w_i
<_{\text{Lex}} v_i$.
\endproclaim

\demo{Proof}  As in [5, Lemma 8], one builds up $p-DF$ and
$<_{\text{Lex}}$ by induction; a term\newline
$pa_0 a_1 \cdots a_{n-1} \circ
a_n$ is put in the set $p-DF$ (and its lexicographic comparison with
terms previously put in is established) only after all the 
iterates $I_m (pa_0 \cdots a_{n-1} ,a_n )$, $m \geq 1$ have been put
in the set.
\enddemo

{\bf Remarks.}  The members of $p-DF$ are terms, and $p-DF$ is closed
under subterms (for $w \leq_L p$, $w$ is the only subterm of $w$, and
for $w = pa_0 a_1 \cdots a_{n-1} \ast a_n$, the subterms of $w$ are $w$
and the subterms of $pa_0 \cdots a_{n-1} ,a_n$).  We will associate
these terms without comment with the members of ${\Cal P}$ they stand
for, when no confusion should arise.  If $w \in p-DF$ and $u$ is a
proper subterm of $w$, then $u <_{\text{Lex}} w$.  Terms of the form
$(u\circ v)w$ or $(u\circ v)\circ w$ are never in $p-DF$.  When using
phrases such as ``$uv \in p-DF$,'' ``$u\circ v \in p-DF$,'' it is
assumed that $u = pa_0 \cdots a_{n-1}$, $v = a_n$ are as in the
definition of $p-DF$ --- isolated exceptions where $uv$ or $u\circ v$
are $\leq_L p$ and are to be considered as singleton terms, will be
noted.

If $u\circ v \in p-DF$, then $u\circ v$ is the
$<_{\text{Lex}}$-supremum of $\{ I_n (u,v) : n \geq 1\}$.

\proclaim{Lemma 3}  The transitivization of the relation 
$\{  \langle u,v\rangle :
u, v \in p -DF$ and either $u$ a proper subterm of $v$, or $v = a\circ
b$ and $u$ is an $I_k (a,b) \}$ is a
well-founded partial ordering $\prec^p$ of $p-DF$.
\endproclaim

\demo{Proof}  Otherwise there would be a sequence $\langle u_n : n <
\omega \}$  with, for each $n$, either $u_{n+1}$ a proper subterm of
$u_n$, or $u_{n+1}$ an iterate of $\langle a,b \rangle$ with $u_n =
a\circ b$, such that no proper subterm of $u_0$ begins such a sequence.
Then $u_0 = r\circ s$, $u_1$ is an iterate of $\langle r,s
\rangle$, and by the nature of such iterates, some $u_n$ must be a
subterm of $r$ or of $s$, a contradiction.
\enddemo

\proclaim{Lemma 4}
\item {(i)}
If $w, a, b_0 ,b_1 ,\ldots ,b_n \in {\Cal P}$, $w b_0 b_1 \cdots b_{n-1}
\ast b_n$ is prenormal with respect to $<_L$, and $b_0 <_L a$, then $wb_0
b_1 \cdots b_{n-1} \ast b_n <_L wa$.
\item {(ii)} For $p \in {\Cal P}$, $u, v \in p-DF$, $u <_{\text{Lex}}
v$ iff $u <_L v$.
\endproclaim

\demo{Proof}  
(i) By induction on $i$ we show $wa >_L wb_0 b_1 \cdots
b_{i-1} \circ b_i$.  For $i=0$, it is Theorem 1(i).  For $i = k+1$, $wa
= (wb_0 \cdots b_{k-1} \circ b_k )u_0 \cdots u_{m-1} \ast u_m \geq_L
wb_0 \cdots b_{k-1} (b_k u_0 ) = wb_0 \cdots b_{k-1} b_k (wb_0 \cdots
b_{k-1} u_0 ) \geq_L wb_0 \cdots b_{k-1} b_k$ ($b_{k+1} r$) for some $r$
(since\newline
$b_{k+1} \leq_L  wb_0 \cdots b_{k-1}$) $>_L wb_0 \cdots b_k \circ
b_{k+1}$.

(ii) It suffices to show $u <_{\text{Lex}} v \Rightarrow u <_L
v$ (the other direction following from that, the linearity of
$<_{\text{Lex}}$ and the irreflexivity of $<_L$).  By induction on
ordinals $\alpha$, suppose it has been proved for all pairs $\langle u'
,v ' \rangle$, $u' ,v' \in p-DF$, such that $u'$ and $v'$  have rank less 
than $\alpha$ with respect to $\prec^p$.  If either of $u, v$ is $\leq_L p$,
or if the associated sequence of $u$ is a proper initial segment of the
associated sequence of $v$, the result is clear.  So, passing to a
truncation $p, a_0 ,a_1 ,\ldots ,a_n$ of $u$'s associated sequence if
necessary, we have $u \geq_{\text{Lex}} pa_0 a_1 \cdots a_n$, $v = pa_0
a_1 \cdots a_{n-1} v_n v_{n+1} \cdots v_{m-1} \ast v_m$, some $m \geq
n$, with $v_n <_{\text{Lex}} a_n$ (the reason why $v$ cannot be $pa_0
\cdots a_{i-1} \circ a_i$ for some $i < n$ is that $a_n
\leq_{\text{Lex}} v_n$ would then hold).  Thus $u \geq_L pa_0 a_1 \cdots
a_n$ (clear), $v_n  <_L a_n$ (by the induction hypothesis) and for each
$i$, $v_{i+1} \leq_L pa_0 a_1 \cdots a_{n-1} v_n \cdots v_{i-1}$ (by the
induction hypothesis).  Then apply part (i) of this lemma.

Thus,for $p,q \in {\Cal P}$, 
to determine which of $q <_L p$, $q = p$, $p
<_L q$ holds, lexicographically compare $|q|^x$ and $|p|^x$.

Write $<_L$ for $<_{\text{Lex}}$ below.  ``Prenormal,'' below, will be
with respect to $<_L$.  For $q, p \in {\Cal P}$, let $|q|^p$ be the
$p-DF$ representation of $q$, if it exists.

Recall that the main theorem is that for all $q, p \in {\Cal P}$ $|q|^p$
exists.  From Lemma 4, this may be stated as a type of division
algorithm:  ~if $q, p \in {\Cal P}$ and $p <_L q$, then there is $a
<_L$-greatest $a_0 \in {\Cal P}$ with $pa_0 \leq_L q$, and if $pa_0 <_L
q$, then there is $a <_L$-greatest $a_1 \in {\Cal P}$ with $pa_0 a_1
\leq_L q$, etc., and for some $n$, $pa_0 a_1 \cdots a_n = q$ or $pa_0
a_1 \cdots a_{n-1} \circ a_n = q$.  And, if this process is repeated
for each $a_i$, getting either $a_i \leq_L p$ or $a_i = pa_i^0  a_i^1
\cdots a_i^{m-1} \ast a_i^m$, and then for each $a_i^k$, etc., then the
resulting tree is finite.  The normal form theorems in [5] correspond to
similar algorithms --- they were proved there just for $p \in {\Cal A}$,
and the present form has their generalizations to all $p \in {\Cal P}$
as a corollary.

In certain cases on $u, v \in p-DF$ (when ``$u \sqsupset^p  v$''), the
existence of $|uv|^p , |u\circ v|^p$ can be proved directly.  We define
$u \sqsupset^p v$ by induction:  ~suppose $u' \sqsupset^p w$ has been
defined for all proper subterms $u'$ of $u$ and all $w \in p-DF$.
\item {(i)} If $u <_L p$, then $u \sqsupset^p v$ iff $u >_L v$ and
$u\circ v \leq_L p$.
\item {(ii)} $p \sqsupset^p v$ for all $v$.
\item {(iii)} $pa \sqsupset^p v$ iff $v \leq_L p$ or $v = pa_0 a_1
\cdots a_{n-1} \ast a_n$ with $a \sqsupset^p a_0$; $p\circ a
\sqsupset^p v$ for all $v$.
\item {(iv)} For $n \geq 1$, $pa_0 a_1 \cdots a_n \sqsupset^p v$ iff
either $v \leq_L pa_0 a_1 \cdots a_{n-1}$ or\newline
$v = pa_0 a_1 \cdots
a_{n-1} ~v_n v_{n+1} \cdots v_{i-1} \ast v_i$ with $a_n \sqsupset^p v_n$
and $a_n \circ v_n \leq_L pa_0 a_1 \cdots a_{n-2}$.
\item {(v)} For $n \geq 1$, $pa_0 a_1 \cdots a_{n-1} \circ a_n
\sqsupset^p v$ iff $a_n \sqsupset^p v$ and $a_n \circ v \leq_L pa_0 a_1
\cdots a_{n-2}$.

\proclaim{Lemma 5}  If $u \sqsupset^p v$, $w \in p-DF$, and $v \geq_L
w$, then $u \sqsupset^p w$.
\endproclaim

\demo{Proof}  By induction on $u$ in $p-DF$.
\enddemo

\proclaim{Lemma 6}  If  $u \sqsupset^p v$, then $|uv|^p$ and $|u\circ
v|^p$ exist, and $|uv|^p \sqsupset^p u$.
\endproclaim

\demo{Proof}  Assume the lemma has been proved for all $\langle u'
,w\rangle$, $w \in p-DF$ and $u'$ a proper subterm of $u$, and for all
$\langle u, v' \rangle$, $v'$ a proper subterm of $v$.  Suppose $u
\sqsupset^p v$.
\item {(i)} $u <_L p$.  Then $uv <_L u\circ v \leq_L p$, so $uv$ and
$u\circ v$, as terms of length one, are in $p-DF$, and $uv\circ u =
u\circ v$, so similarly $uv \sqsupset^p u$.
\item {(ii)} $u = p$.  Then $|pv|^p = pv$, $|p\circ v|^p = p\circ v$,
and $pv \sqsupset^p p$.
\item {(iii)} $u = pa$.  Then if $v \leq_L p$ it is clear, so assume 
$v = pb_0 b_1 \cdots b_{n-1} \ast b_n$, where $a \sqsupset^p b_0$.  The 
cases are 
\itemitem {(a)} $v = pb$.  Then $|uv|^p = p |ab|^p$, $|u\circ v|^p = p
|a\circ b|^p$, when $|ab|^p$, $|a\circ b|^p$ exist by induction.  And
since by induction $|ab|^p \sqsupset a$, we have $|uv|^p \sqsupset u$.
\itemitem {(b)} $v = p\circ b$.  Then $|uv|^p = |pa(p\circ b)|^p = p
|a\circ b|^p p \circ p|ab|^p$ by the induction hypothesis and Theorem
1(i).  Similarly $|u\circ v|^p = |pa\circ (p\circ b)|^p = p\circ |a\circ
b|^p$.  To see $|uv|^p \sqsupset^p u$, we have $p |ab|^p \sqsupset_p
pa$, as $|ab|^p \sqsupset^p a$ holds by the induction hypothesis, and
$p(ab)\circ pa = p(ab\circ a) = p(a\circ b)$.
\item {(c)} $v = pb_0 b_1 \cdots b_{n-1} \ast b_n$ for $n \geq 1$.  Then
$|uv|^p = p|a\circ b_0 |^p b_1 |pab_2 |^p \cdots |pab_{n-1} |^p \ast
|pab_n |^p$ by the induction hypothesis and Theorem 1(i) and Lemma
4(ii).  And in the case $\ast = \cdot$, $|u\circ  v|^p = |uv\circ u|^p
= p|a\circ b_0 |^p b_1 |pab_2 |^p \cdots |pab_n |^p \circ pa$.  In the
case $\ast = \circ$, $|u\circ v|^p = |uv\circ u|^p = p|a\circ b_0 |^p b_1
|pab_2 |^p \cdots |pab_{n-1}  |^p \circ |pab_n \circ pa|^p$, namely,
$|pab_n \circ pa|^p = |pa\circ b_n |^p$ exists by induction and is $<_L
p(a\circ b_0 )b_1 (pab_2 )\cdots (pab_{n-2} )$ by $b_n <_L pb_0 \cdots
b_{n-2}$ and Theorem 1(i).  To see $|uv|^p \sqsupset^p u$, it is immediate
if $\ast = \cdot$, and if $\ast = \circ$, $pab_n \sqsupset^p pa$ by 
induction, and $pab_n \circ pa = pa\circ b_n <_L pa(pb_0 \cdots b_{n-2} ) 
= p(a\circ b_0 )b_1 (pab_2 )\cdots (pab_{n-2} )$, as desired.
\item {(iv)} $u = pa_0 a_1 \cdots a_n$, $n \geq 1$.  Then the case where
the induction hypothesis is used is where $v = pa_0 a_1 \cdots
a_{n-1} b_n \cdots b_{m-1} \ast b_m$, where $a_n \sqsupset^p b_n$ and
$a_n \circ b_n \leq_L pa_0 a_1 \cdots a_{n-2}$.  The cases and
computations are similar to (iii).
\item {(v)} $u = pa_0 a_1 \cdots a_{n-1} \circ a_n$, $n \geq 1$.  Then
$a_n \sqsupset^p v$, so $|a_n v|^p$, $|a_n \circ v|^p$ exist, and $a_n v
<_L a_n \circ v \leq_L pa_0 \cdots a_{n-2}$.  Thus $|uv|^p = pa_0 a_1
\cdots a_{n-1} |a_n v|^p$, $|u\circ v|^p = |pa_0 a_1 \cdots a_{n-1}
\circ |a_n \circ v|^p |^p$ which is $pa_0 \cdots a_{i-1} \circ |a_i
\circ v|^p$ , where $i \leq n$ is greatest such that $i = 1$ or $a_i
\circ v <_L pa_0 \cdots a_{i-1}$.
And for $|uv|^p \sqsupset^p u$, we have $|a_n v|^p
\sqsupset^p a_n$, and $a_n v \circ a_n = a_n \circ v \leq_L pa_0 \cdots
a_{n-2}$, as desired.

\proclaim{Lemma 7}  Suppose $p, q \in {\Cal P}$, $w \in q-DF$.  Then
\item {(i)} $|pw|^{pq}$ exists, and $|pw|^{pq} \sqsupset^{pq} p$.
\item {(ii)} If $|pw|^{p\circ q}$ exists, then $|pw|^{p\circ q}
\sqsupset^{p\circ q} p$.
\endproclaim

\demo{Proof}  We check part (ii), part (i) being similar.  Assume the
lemma is true for all proper components $w'$ of $w$.  If $w \leq q$,
then $pw < p\circ w \leq p \circ q$ and, by $pw \circ p = p\circ w$, we
have $pw \sqsupset^{p\circ q} p$.  So assume the most general case on
$w$, $w = qa_0 a_1 \cdots a_{n-1} \circ a_n$.  Then $pw = (p\circ q)a_0
(pa_1 )\cdots (pa_{n-1} )\circ (pa_n )$ is prenormal, so if
$|pw|^{p\circ q}$ exists, then by Lemma 4(i) and (ii) $|pw|^{p\circ q} 
= (p\circ q)|a_0 |^{p\circ q} |pa_1 |^{p\circ q} \cdots |pa_{n-1}
|^{p\circ q} \circ |pa_n |^{p\circ q}$.  Then $|pa_n |^{p\circ q}
\sqsupset^{p\circ q} p$ by the induction assumption, and $pa_n \circ p
= p\circ a_n <_L p(qa_0 \cdots a_{n-2} ) = (p \circ q) a_0 (pa_1 )\cdots
(pa_{n-2} )$.  So $|pw|^{p\circ q} \sqsupset^{p\circ q} p$.  The case
$n=0$ yields $p(q\circ a) = p(qa\circ q) = (p\circ q)a \circ (pq)$ and
is similarly checked, using that $pq \sqsupset^{p\circ q} p$.
\enddemo

Note, for $F$ a finite subset of ${\Cal P}$, the following induction
principle:  ~if $S \subseteq {\Cal P}$, $S \not= \emptyset$, then there
is a $w \in S$ such that for all $u$, if $pu \leq w$ for some $p \in F$,
then $u \notin S$.  Otherwise some $w \in S$ would be $\geq$
arbitrarily long compositions of the form $p_0 \circ p_1 \circ \cdots
\circ p_n$, each $p_i \in F$.  By Theorem 1(ii), some $p \in F$ would 
occur at least $m$ times in one of
these compositions, where $p^m > p^{(m)} > w$, and applications of
the $a \circ b = ab \circ a$ law would give $p^m \leq p_0 \circ \cdots
\circ p_n \leq w$, a contradiction to Theorem 1(v) and (viii).

\proclaim{Theorem}  For all $w, r \in {\Cal P}$, $|w|^r$ exists.
\endproclaim

\demo{Proof}  We show that $T = \{ r \in {\Cal P} : \text{for all} ~w
\in {\Cal P} , ~|w|^r ~\text{exists} \}$ contains $x$ and is closed under
$\cdot$ and $\circ$.
\item {(i)} $x \in T$.  Suppose, letting $F = \{ x\}$ in the induction
principle, that $|w|^x$ does not exist but $|u|^x$ exists for all $u$ such
that $xu \leq w$.  Pick $v \leq w$ such that $|v|^x$ doesn't exist, and,
subject to that, the $(x,x)$-normal form of $v$ ([5], Lemmas 25, 27,
Theorem 28) has minimal length.  The $(x,x)$-normal form of $v$ is a
term $xa_0 a_1 \cdots a_{n-1} \ast a_n$, which is prenormal, where $a_0$
is in the normal form of [5] (see the corollary below), and for $i >
0$, each $a_i$ is in $(x,x)$-normal form.  Then for $i > 0$, each $|a_i
|^x$ exists, and since $xa_0 \leq w$, $|a_0 |^x$ exists.  Thus $|v|^x$
exists, $|v|^x = x|a_0 |^x \cdots |a_{n-1} |^x \ast |a_n |^x$.
\item {(ii)} $p, q \in T$ implies $pq \in T$.  For $u \in p-DF$, define
the $\langle p,q\rangle -DF$ of $u$ as follows.  If $u \leq p$, the
$\langle p,q\rangle -DF$ of $u$ is $u$.  If $u = pa_0 a_1 \cdots a_{n-1}
\ast a_n$, the $\langle p,q\rangle -DF$ of $u$ is $p\bar a_0 \bar a_1
\cdots \bar a_{n-1} \ast \bar a_n$, where $\bar a_0 = |a_0 |^q$ and for
$i > 0$, $\bar a_i$ is the $\langle p,q\rangle -DF$ of $a_i$.  Then by
assumption every $r \in {\Cal P}$ has a $\langle p,q\rangle -DF$
representation.  Pick $v$ such that $|v|^{pq}$ doesn't exist, and
subject to that, the $\langle p,q\rangle -DF$ representation of $v$ has
minimal length.  If $v \leq pq$, we are done.  So assume $v$'s $\langle
p,q\rangle -DF$ representation is $p(qa_0 a_1 \cdots a_{n-1} \ast a_n )
b_0 b_1 \cdots b_{m-1} \ast b_m$, where the proof for $n \geq 0$ and the
first $\ast$ being $\circ$ will cover all cases.  Then $v = pq(pa_0
)(pa_1 )\cdots (pa_{n-1} )(pa_n b_0 ) b_1 \cdots b_{m-1} \ast b_m$.
Then $|pa_0 |^{pq} \cdots |pa_{n-1} |^{pq} , |pa_n |^{pq} , |b_0 |^{pq}
\cdots |b_m |^{pq}$ all exist by the minimality of $v$'s $\langle p,q
\rangle -DF$ representation.  And since $b_0 \leq p$, $|pa_n |^{pq}
\sqsupset^{pq} b_0$ by Lemma 7(i), and $|pa_n b_0 |^{pq}$ exists by
Lemma 6.  The sequence $(pq), (pa_0 ) \cdots (pa_{n-1} ), (pa_n b_0 ),
b_1 \cdots b_{n-1} ,b_n$ need not be prenormal.  But we claim $|p(qa_0
\cdots a_{n-1} \circ a_n )|^{pq} = pq |pa_0 |^{pq} \cdots |pa_{n-1}
|^{pq} \circ |pa_n |^{pq} \sqsupset^{pq} |b_0 |^{pq}$.  The equality is
clear.  For the $\sqsupset^{pq}$ relation, we have $|pa_n |^{pq}
\sqsupset^{pq} |b_0 |^{pq}$ and $pa_n \circ b_0 \leq pa_n \circ p = p
\circ a_n \leq pq (pa_0 ) \cdots (pa_{n-2} )$ since $a_n < pa_0 \cdots
a_{n-2}$, giving the claim.  So by Lemma 6, $|p(qa_0 \cdots a_{n-1}
\circ a_n ) b_0 |^{pq} \sqsupset^{pq} |p(qa_0 \cdots a_{n-1} \circ a_n
)|^{pq} \geq b_1$.  By Lemma 5, $|p(qa_0 \cdots a_{n-1} \circ a_n ) b_0
|^{pq} \sqsupset^{pq} |b_1 |^{pq}$.  With this as the first step,
iterate Lemma 6 and Lemma 5, $m$ times, to get that $|p(qa_0 \cdots
a_{n-1} \circ a_n ) b_0 b_1 \cdots b_{m-1} \ast b_m |^{pq}$ exists.
\item {(iii)} $p, q \in T$ implies $p \circ q \in T$.  Letting $F = \{
q\}$ in the induction principle, suppose $|w|^{p\circ q}$ does not exist
but $|a|^{p\circ q}$ exists for all $a$ such that $qa \leq w$.  Pick $v
\leq w$ such that $|v|^{p\circ q}$ does not exist and, subject to that,
the $\langle p, q\rangle -DF$ representation of $v$ has minimal length.
If $v \leq p\circ q$, then again the cases on the $\langle p,q\rangle -DF$
representation of $v$ are covered by the proof where that representation
is $p(qa_0 \cdots a_{n-1} \circ a_n )b_0 b_1 \cdots b_{m-1} \ast b_m$.

Then $v = (p\circ q)a_0 (pa_1 )\cdots (pa_{n-1} )(pa_n b_0 )b_1 \cdots
b_{m-1} \ast b_m$.  As in case (ii),\newline
$|pa_1 |^{p\circ q} ,\ldots
,|pa_{n-1} |^{p\circ q} ,|pa_n |^{p\circ q} ,|b_0 |^{p\circ q} ,\ldots
,|b_m |^{p\circ q}$ exist, and using Lemma 7(ii) and\newline
Lemma 6, $|pa_n b_0
|^{p\circ q}$ exists.  And since $qa_0 \leq v$, $|a_0 |^{p\circ q}$
exists by the induction principle.  Thus $|p(qa_0 \cdots a_{n-1} \circ
a_n )|^{p\circ q}$ exists, and, as in case (ii), is $\sqsupset^{p\circ
q} |b_0 |^{p\circ q}$.  Then iterate Lemmas 6 and 5 as in case (ii) to
obtain the existence of $|v|^{p\circ q}$.  This completes the proof of
the theorem.

For $p \in {\Cal P}$, say that a term $w$ in the alphabet $\{ q : q <_L
p\} \cup \{ p^{(i)} : i < \omega \}$ is in $p$-normal form $(p-NF)$ if
either $w <_L p$ is a term of length one, or $w = p^{(i)} a_0 a_1 \cdots
a_{n-1} \ast a_n$, where each $a_k \in p-NF$, $p^{(i)} a_0 a_1 \cdots
a_{n-1} \ast a_n$ is prenormal, and $a_0 <_L p^{(i)}$.  Let $|w|_p$ be
the $p-NF$ representation of $w$ if it exists.  As in [5], Lemmas 9 and
12, such a representation is unique.  It is proved in [5] that for all
$p \in {\Cal A}$ and $w \in {\Cal P}$, $|w|_p$ exists.  The $DF$ theorem
allows this to be extended to $p \in {\Cal P}$.

\proclaim{Corollary}  If $p, w \in {\Cal P}$, then $|w|_p$ exists.
\endproclaim

\demo{Proof}  By induction on $w \in p-DF$.  If $w <_L p$, we are done;
so assume $w$ is the $p-DF$ term $pa_0 a_1 \cdots a_{n-1} \ast a_n$.
Then each $|a_i |_p$ exists, and if $a_0 <_L p$, we are done.  Also, if
$a_0 = p$, then the $p-NF$ expression for $w$ is $p^{(1)} |a_1 |_p
\cdots |a_{n-1} |_p \ast |a_n |_p$.  Without loss of generality assume
$a_0$'s $p-NF$
representation is $p^{(m)} b_0 b_1 \cdots b_{k-1} \circ b_k$.  Then it
is easily checked that $|pa_0 |_p = p^{(m+1)} |pb_0 |_p \cdots |pb_{k-1}
|_p \circ |pb_k |_p$.  Thus $w = p^{(m+1)} (pb_0 ) \cdots (pb_{k-1} )(pb_k
a_1 )a_2 \cdots a_{n-1} \ast a_n$.  In [4, Theorem 16], a $\sqsupset_p$
theorem is proved for $p-NF$ (for $p \in {\Cal A}$, but a similar result
holds for all $p \in {\Cal P}$).  We may use a version of it, and an 
analog of Lemma 7
above, as Lemmas 6 and 7 were used in Theorem 8, to obtain $|pa_0 |_p
\sqsupset_p a_1$, and then iterate to get the existence of $|w|_p$.  The
details are left to the reader.
\vskip 2em
\baselineskip=12pt

\centerline {\bf References}
\vskip 1em
\item {[1]} P. Dehornoy, Free distributive groupoids, {\it Journal of
Pure and Applied Algebra} {\bf 61} (1989), 123--146.
\vskip 10pt
\item {[2]} P. Dehornoy, Sur la structure des gerbes libres, C.R.A.S.
Paris, t.308, S\'erie I (1989), 143--148.
\vskip 10pt
\item {[3]} P. Dehornoy, Braid groups and left distributive operations,
preprint.                          
\vskip 10pt
\item {[4]} R. Dougherty, On critical points of elementary embeddings,
Handwritten notes, 1988.
\vskip 10pt
\item {[5]} R. Laver, On the left distributive law and the freeness of
an algebra of elementary embeddings, {\it Advances in Mathematics}
{\bf 91} (1992), 209-231.
\vskip 10pt  
\item {[6]} R. Laver, On an algebra of elementary embeddings,preprint.
\vskip 10pt
\item {[7]} J. Steel, On the well foundedness of the Mitchell order,preprint.

\end